\documentclass[12pt,a4]{article}

\usepackage[all]{xy}
\usepackage{dsfont}

\usepackage{amsthm,amsfonts}

\usepackage{graphicx}

\xyoption{all}
\newtheorem{teo}{Theorem}[section]
\newtheorem{prop}[teo]{Proposition}

\newtheorem{cor}[teo]{Corollary}
\newtheorem{lem}[teo]{Lemma}
\newtheorem{rem}[teo]{Remark}

\newcommand{\cO}{\mathcal{O}}
\newcommand{\cM}{\mathcal{M}}
\newcommand{\ra}{\rightarrow}
\newcommand{\CC}{\mathbb{C}}
\newcommand{\PP}{\mathbb{P}}
\newcommand{\LL}{\lambda}
\begin{document}

\hoffset = -1.2truecm

\voffset = -1.8truecm

\title{\bf{Moduli spaces of quadratic complexes and their singular surfaces}}

\author{\textbf{ D. Avritzer}\thanks{Both authors would like
to thank the DAAD (Germany) and CAPES (Brasil) for support during
the preparation of this paper-Project \#191/04} \and \textbf{  H.
Lange}$^{*}$}

\date{}

\maketitle
\begin{abstract}
\noindent
We construct the coarse moduli space $\cM_{qc}(\sigma)$ of quadratic line complexes with a fixed Segre symbol $\sigma$ 
as well as the moduli space $\cM_{ss}(\sigma)$ of the corresponding singular surfaces. We show that the map associating to a quadratic 
line complex its singular surface induces a morphism $\pi: \cM_{qc}(\sigma) \ra \cM_{ss}(\sigma)$. Finally we deduce that the varieties 
of cosingular quadratic line complexes are almost always curves.   
\end{abstract}

\section{Introduction}

Let $G$ denote a smooth quadric in $\PP^5$ over the field of complex numbers, considered as the Pl\"ucker 
quadric parametrizing lines in $\PP^3$. A {\it quadratic complex} or to be more precise a {\it quadratic line complex} 
is by definition a complete intersection $X = F \cap G$ with a quadric $F \subset \PP^5$ different from $G$. 
A quadratic complex can be considered as a pencil of quadrics. Hence the Segre symbol $\sigma = \sigma(X)$ 
of a quadratic complex is well defined. For the definition see Section 2. The first aim of the paper is to construct the moduli spaces of quadratic 
complexes with a fixed Segre symbol.

Let $SO(G)$ denote the special orthogonal group associated to the quadric $G$. Two quadratic complexes $X_1$ and $X_2$ 
are isomorphic if and only if there is a matrix $A \in SO(G)$ such that $X_2 = A(X_1)$. This gives an action of $SO(G)$ 
on the space of quadratic complexes and the notion of semistable quadratic complexes is well defined. It turns out 
(Corollary \ref{cor3.3}) that a quadratic complex $X$ is semistable if and only its discriminant admits at 
least two different roots or equivalently if its Segre symbol consists of at least two brackets. On the other hand, 
a quadratic complex is non-reduced (respectively reducible) if and only if its Segre symbol contains a bracket of length 5 
(respectively 4). Hence all irreducible and reduced quadratic complexes are semistable. The irreducible and reduced 
quadratic complexes with a fixed Segre symbol $\sigma$ form an irreducible variety on which the group $SO(G)$ acts 
and such that all orbits are of the same dimension. This implies that the coarse moduli space $\cM_{qc}(\sigma)$ 
of these quadratic complexes exists and is a quasiprojective variety of dimenson $r-2$, where $r$ denotes the 
number of brackets of the Segre symbol $\sigma$ (see Theorem \ref{teo4.7}).\\

The classical method for investigating quadratic complexes is to study the set of lines in $\PP^3$ parametrized by the points of $X$. 
For any point $p \in \PP^3$, the lines in $\PP^3$ passing through $p$ are parametrized by a plane $\alpha(p)$ containd in $G$.
For a general point $p \in \PP^3$ the intersection $\alpha(p) \cap F$ is a smooth conic. The set 
$S =S(X) = \{p \in \PP^3: rk(\alpha(p) \cap F) \leq 2 \}$ is a surface in $\PP^3$, not necessarily irreducible. 
It is called the {\it singular surface in $\PP^3$} associated to the quadratic complex $X$. The surface $S(X)$ is a quartic
whose singularities depend on the Segre symbol $\sigma(X)$. The next aim of the paper is to construct the moduli space of 
these quartics.

The group $SL(4) = SL(4,\CC)$ acts in a natural way on the space of quartics in $\PP^3$. Hence it makes sense to talk 
about semistability of these quartics. We show (see Proposition \ref{prop6.2}) that a quartic surface in $\PP^3$ is 
semistable with respect to the action of $SL(4)$ if and only if it does not admit a triple point whose tangent cone 
is a cone over a cuspidal plane cubic (possibly degenerated). We use this to construct the moduli space $\cM_{ss}(\sigma)$ of quartics 
with a singularity of type $\sigma$ (Theorem \ref{teo6.4}).

The moduli space $\cM_{qc}(\sigma)$ is constructed using the action of the group $SO(G) \simeq SO(6)$, whereas the moduli space
$\cM_{ss}(\sigma)$ is constructed with respect to the action of the group $SL(4)$. Note that there is an isomorphism
$$
PSO(6) \simeq PSL(4).
$$
We use this isomorphism to show that the map associating to a quadratic complex $X$ its singular surface $S(X)$ induces a morphism 
$$
\pi: \cM_{qc}(\sigma) \longrightarrow  \cM_{ss}(\sigma).
$$
In the classical literature (see e.g. \cite{J}) two quadratic complexes are called {\it cosingular} if their singular 
surfaces are isomorphic. So the fibres of the morphism $\pi$ are just the varieties of cosingular complexes. In \cite{K}
Klein showed the variety of cosingular complexes in the generic case, i.e. for $\sigma = [111111]$, is just the projective line.
We show (see table 7.3) that the varieties of cosingular complexes are generically curves, except in one case. Finally, in the last section 
we reprove Klein's result using our set up.\\ 
 
We would like to thank G.-M. Greuel for pointing out the paper \cite{FK} to us.\\ 

{\it Notation}: A quadric $F \in \PP^n$ is up to a nonzero
constant given by a symmetric matrix of degree $n+1$. By a slight
abuse of notation we denote the matrix by the same letter as the
quadric it defines. We work over the field of complex numbers. Most of the results are valid over an arbitrary algebraically 
closed field of characteristic $\neq 2$, however some of the computations use complex numbers. Hence the groups $SO(G), SO(6), SL(4)$
etc. always mean the corresponding complex Lie groups.  \\

\section{The Segre symbol of a quadratic complex}

Let $G \subset \PP^5$ be a smooth quadric which we fix in the sequel.
We consider $G$ as the Pl\"ucker quadric, although not always with
Pl\"ucker coordinates. In fact, we choose the coordinates
appropriate to the statement we want to prove. However in any
case, a point $x \in G$ represents a line in $\PP^3$, denoted by
$l_x$. Let $F$ denote a quadric in $\PP^5$, different from $G$.
The complete intersection
$$
X = F \cap G,
$$
parametrizes a set of lines in $\PP^3$, which is classically
called a quadratic line complex. We call the variety $X$ itself
the {\it quadratic line complex} or to be short the {\it quadratic complex} defined by $F$.

A quadratic line complex determines a pencil of quadrics in
$\PP^5$, namely
$$
{\cal P} = \{Q_{(\lambda:\mu)} = \lambda F + \mu G\; | \;
(\lambda:\mu) \in \PP^1 \}
$$
which we call the pencil {\it associated} to the quadratic line
complex. Note that the space of pencils of quadrics is by
definition the Grassmannian $Gr(2,21)$ of lines in $\PP^{20}$ and
thus of dimension 38, whereas, as we shall see, the space of
quadratic line complexes is of dimension 19 only. Pencils of
quadrics are classified according to their Segre symbol. The {\it
Segre symbol of the quadratic complex} is by definition the Segre
symbol of the associated pencil of quadrics.

Let us recall the definition of the  Segre symbol: The {\it
discriminant of the pencil} ${\cal P}$ is by definition the binary
sextic
$$
 \Delta=\Delta(\lambda,\mu ):=\det(\lambda F + \mu
G).
$$
The discriminant $\Delta$ of ${\cal P}$ depends on the choice of
the matrices $F$ and $G.$ The roots of $\Delta,$ however,  are
uniquely determined up to an isomorphism of $\PP^1$. In particular,
the multiplicities of the roots are uniquely determined.

Suppose $({\bar \lambda}:{\bar \mu})$ is a root of $\Delta.$ It
may also happen that all the subdeterminants of ${\bar \lambda} F
+{\bar \mu} G $ of a certain order vanish. Suppose that all
subdeterminants of order $6-d$ vanish for some $d \geq 0,$ but not
all subdeterminants of order $5-d.$ This means that the quadric
$Q_{({\bar \lambda}:{\bar \mu})}$ is a d-cone with vertex a linear
space of dimension $d$ and directrix a smooth quadric in a linear
subspace of dimension $4-d$ in $\PP^5$.

Let  $l_i$ denote the minimum multiplicity of the root
 $({\bar \lambda}:{\bar \mu})$ in the subdeterminants of order $6-i,$ for
$i=0,1,\dots,d.$ Then $l_i > l_{i+1}$ for all $i$ so that
$e_i:=l_i-l_{i+1} > 0,$ and we have :
$$
 \Delta(\lambda,\mu) =
(\lambda {\bar \mu} - {\bar \lambda} \mu)^{e_0} \dots (\lambda
{\bar \mu} - {\bar \lambda} \mu)^{e_{d}} \Delta_1(\lambda,\mu),
$$ with
$\Delta_1({\bar \lambda},{\bar \mu}) \neq 0.$ The numbers $e_i$
are called the {\it characteristic numbers} of the
 the root $({\bar \lambda}:{\bar \mu})$ and the
 factors $(\lambda {\bar \mu} -{\bar \lambda} \mu)^{e_i}$ are called
the {\it elementary divisors} of the pencil ${\cal P}.$

If $(\lambda_i:\mu_i)\mbox{ for }i=1,\dots,r $ are the  roots  of
$\Delta$ and $e_0^j,\dots, e_{d_j}^j$ the characteristic numbers
associated to the root $(\lambda_j:\mu_j)$ and $d_1 \geq d_2 \geq
\dots \geq d_r,$ then
$$
\sigma_X = \sigma_{\cal P} =
[(e_0^1\dots e_{d_1}^1)(e_0^2\dots e_{d_2}^2)
\dots(e_0^r \dots e_{d_r}^r)]
$$
is called the {\it Segre symbol} of the quadratic complex $X$ or
the pencil ${\cal P}.$ The parentheses are omitted if $d_i=1.$
In order to make it unique, we assume that the expressions
$(e_0^i,\ldots,e_{d_i}^i)$ are ordered lexicographically if $d_i =
d_j$. We call these expressions the {\it brackets} of the Segre
symbol $\sigma_X$ or of the pencil ${\cal P}$ (even if the
parentheses are omitted, i.e. $d_i =1$).

It is a classical fact (see e.g. \cite[p. 278]{HP}) that 2 pencils
of quadrics ${\cal P}_1$ and ${\cal P}_2$  in $\PP^n$, whose
discriminants
 have roots exactly
 at  $(\lambda_i^1:\mu_i^1)$
and $(\lambda_i^2:\mu_i^2),$ are isomorphic, that is, projectively
equivalent in $\PP^n$, if and only if they have the same Segre
symbol and there is an automorphism of $\PP^1$ taking
$(\lambda_i^1:\mu_i^1)$ to  $(\lambda_i^2:\mu_i^2)$ for all i, where the brackets corresponding to 
$(\lambda_i^1:\mu_i^1)$ and $(\lambda_i^2:\mu_i^2)$ are of the same type.
This can be used  to define a normal form for those pencils ${\cal
P}$, whose discriminant is not identically zero (see \cite[p.
280]{HP}):

For every $e^i_j$ occurring in the Segre symbol of $X$ consider
the $e^i_j \times e^i_j$-matrices
$$
F_{ij} =\left(\begin{array}{ccccc}0&0&\dots&1& \frac{\lambda_i}{\mu_i}\\
                                        0&\dots&1& \frac{\lambda_i}{\mu_i}&0\\
                                        \dots&\dots&\dots&\dots&\dots\\
                                        1&\frac{\lambda_i}{\mu_i}&0&\dots&0\\
                                        \frac{\lambda_i}{\mu_i}&0&0&\dots&0
\end{array}\right)\;
\mbox{ and }\;
 G_{ij}=\left(\begin{array}{ccccc}0&0&\dots&0&1\\
                                          0&0&\dots&1&0\\
                                            \dots&\dots&\dots&\dots&\dots\\
                                             0&1&\dots&0&0\\
                                             1&0&\dots&0&0
\end{array}\right).
$$
The coordinates of $\PP^5$ can be chosen in such a way that $F$
and $G$ are given as block diagonal matrices as follows
$$
F = \mbox{diag} (F_{11},\cdots,F_{rd_r}) \quad \mbox{and} \quad G
= \mbox{diag} (G_{11},\cdots,G_{rd_r}).
$$
We call these coordinates {\it Segre coordinates} of the quadratic
complex $X$. Note that Segre coordinates are not uniquely
determined.

\begin{rem}
{\em If $F$ and $G$ are given in Segre coordinates, the matrix
$(G^{-1}F)^t$ will be in Jordan normal form. This gives another way to
determine the Segre normal form: If $X = F \cap G$, choose the
coordinates in such a way that the matrix $(G^{-1}F)^t$ is in Jordan
normal form. Then the Segre normal form can be read off from
this.}
\end{rem}

From the Segre normal form it is easy to derive the following
lemma.

\begin{lem}
Let $(e_0,\cdots,e_d)$ denote a bracket in the Segre symbol of a
quadratic complex. Then $e_0 \geq e_1 \geq \cdots \geq e_d$.
\end{lem}

The lemma is valid for any pencil of quadrics in $\PP^n$, a
general quadric of which is smooth. The proof is the same.

\begin{proof}
The Segre symbol of a pencil $\lambda F + \mu G$ does not depend
on the coordinates chosen. Hence we may choose Segre coordinates
for the pencil. In particular, we can choose the coordinates in
such a way, that the matrix ${\bar F}$ of the block in the matrix
$F$ corresponding to the root $({\bar \lambda}:1)$, which belongs
to the bracket $(e_0, \cdots,e_d)$, is of the form
$$
{\bar F} = \mbox{diag}({\bar F}_0, \cdots, {\bar F}_d) \quad
\mbox{with} \quad {\bar F}_i = \left(\begin{array}{ccccc}0&0&\dots&1& {\bar \lambda}\\
                                        0&\dots&1& {\bar \lambda}&0\\
                                        \dots&\dots&\dots&\dots&\dots\\
                                        1&{\bar \lambda}&0&\dots&0\\
                                        {\bar \lambda}&0&0&\dots&0
\end{array}\right)
$$
such that the sizes ${\bar e}_i$ of the matrices ${\bar F}_i$ are
ordered as follows: ${\bar e}_0 \geq {\bar e}_1 \geq \cdots \geq
{\bar e}_d$. We have to show that $e_i = {\bar e_i}$ for all $i$.
For this it suffices to show that the minimum multiplicity of the
root ${\bar \lambda}$ in the subdeterminants of order $6-i$ of the
matrix ${\bar F}$ is $l_i= \sum_{j=i}^d {\bar e_j}$. Clearly we
have $l_i \geq \sum_{j=i}^d {\bar e_j}$ for all $i$ and have to
show only that the minimum can be obtained. For $i=0$ this is
clear. For $i=1$ cancel the last line and column of the matrix
$\bar F$ and compute the corresponding minor to see this. Then
proceed successively always cancelling the line and column given
by the last line and column of the submatrix ${\bar F}_i$. The
corresponding minor always is the minimum.
\end{proof}

The quadric line complex $X$ is by definition the base locus of
the pencil ${\cal P}.$ Thus $X$ is the intersection of any two
different quadrics of the pencil.

The brackets in the Segre Symbol correspond 1-1 to the cones in
the pencil. We call the cone $Q_{(\lambda_i:\mu_i)}$ corresponding
to the bracket $(e_0^i,\dots,e_{d_i}^i)$ a cone {\it of type}
$(e_0^i,\dots,e_{d_i}^i).$ The quadric $Q_{(\lambda_i:\mu_i)}$ is
then a $d_i$-cone and the corresponding root in the discriminant
$\Delta$ is a root of multiplicity $e^i:=\sum_{j=0}^{d_i} e^i_j.$
By a slight abuse of notation we call $Q_{(\lambda_i:\mu_i)}$ a
{\it d-cone of multiplicity} $e^i.$ The following table gives a
list of the singularities of $X$ corresponding to the brackets occurring
in this paper.
\begin{center}
\begin{tabular}{c|c|c|c}
bracket& dim of vertex & vertex $ \cap \; X$ & type of singularities\\
\hline \hline 1 & 0 & $ \emptyset$ & no singularities\\
\hline 2 & 0 & 1 point & $A_1$ \\
\hline 3 & 0 & 1 point & $A_2$ \\
\hline 4 & 0 & 1 point & $A_3$ \\
\hline $(11)$ & 1 & 2 different points & $A_1$ \\
\hline $(21)$ & 1 & 1 point & $A_2$ \\
\hline $(22)$ & 1 & 1 point & $A_3$ \\
\hline $(111)$ & 2 & smooth conic $C$ & $X$ singular along $C$\\
\hline $(211)$ & 2 & rank 2 conic $C$ & $X$ singular along $C$ \\
\hline
\end{tabular}
\end{center}

\section{Semistable quadratic complexes}

Recall that we fixed a smooth quadric $G$ in $\PP^5$. Considering
the projective space $\PP^{20}$ as the space parametrizing nontrivial
quadrics in $\PP^5$, a quadratic line complex is given by a line
in $\PP^{20}$ passing through $G$. Thus the space of quadratic
line complexes can be considered as the closed subvariety
$$
LC = \{ L \in Gr(2,21) \;|\; G \in L \}
$$
of the Grassmannian of lines in $\PP^{20}$.

Two quadratic line complexes $X_1$ and $X_2$ are called {\it
isomorphic} if there is an automorphism $A$ of $\PP^5$ with
$X_2 = A(X_1)$ which fixes the quadric $G$. The group of
automorphisms of $\PP^5$ fixing the quadric $G$ is by definition
the group $PSO(G) \simeq PSO(6,k)$. We work instead with the
finite covering $SO(G)$. Hence we get an action of the reductive
group $SO(G)$ on the projective variety $LC$. Since $LC$ clearly
admits an $SO(G)$-linearized line bundle, the notion of
semistability is well defined for quadratic line complexes.

In order to determine the semistable quadratic line complexes, we
will use various coordinates of $\PP^5$. Let us recall the
relation between the corresponding special orthogonal groups. Let
$G$ and $G'$ denote the matrices of the Pl\"ucker quadric with
respect to two coordinate systems. We normalize the matrices such
that the determinants of $G$ and $G'$ are $1$ and denote
the corresponding groups by $SO(G)$ and $SO(G')$. Let $A$ denote a
matrix of the coordinate change, so that $A^t G A = G'$ and $A^t =
A^{-1}$. Then
\begin{equation} \label{eq2.1}
SO(G) \ra SO(G'), \qquad g \mapsto
h=A^t g A
\end{equation}
is an isomorphism of groups.

Now choose the coordinates in such a way that $G$ is given by the
matrix $1_6$. In classical terminology the corresponding
coordinates are called {\it Klein coordinates}. Here the
corresponding group $SO(G)$ coincides with the usual orthogonal
group $SO(6)$. Let $S_0$ denote the space of quadrics in
$\PP^5$ with trace 0, i.e. the vector space of nonzero symmetric
$6 \times 6$-matrices of trace 0 modulo $\CC^*$. Obviously $S_0
\simeq \PP^{19}$ and the group $SO(6)$ acts on $S_0$ by $(g,M)
\mapsto g^tMg$.

\begin{prop} \label{prop3.1}
There is a canonical isomorphism
$$
\Phi: LC \ra S_0
$$
compatible with the actions of $SO(G)$ and $SO(6)$. In particular
the variety of quadratic complexes is isomorphic to $\PP^{19}$.
\end{prop}

Note that that the coordinates for $LC$ can be chosen arbitrarily.
That is the reason for denoting the group acting on $LC$ by
$SO(G)$.

\begin{proof}
Given any quadratic complex $X$, choose Klein coordinates. Then
the associated pencil of quadrics $\{\lambda F + \mu G \;|\;
(\lambda:\mu) \in \PP^1 \}$ contains exactly one quadric of trace
0, namely
$$
F_0 = F - \frac{tr F}{6}G.
$$
Certainly this definition does not depend on the choice of $F$.
Conversely, if $F_0$ is a non-zero quadratic form of trace 0, then
$F_0$ and $G \;(=1_6)$ are linearly independent and thus determine a
quadratic line complex. Certainly the maps $X \mapsto F_0$ and
$F_0 \mapsto X$ are algebraic and inverse to each other, giving
the canonical isomorphism as stated.

In both cases the special orthogonal group acts by conjugation and
according to (\ref{eq2.1}) this is independent of the chosen
coordinates. Hence the maps are compatible with the given actions.
\end{proof}
As a consequence of Proposition \ref{prop3.1} we get that a
quadratic complex is semistable with respect to the action of
$SO(G)$ if and only if the associated quadric of trace 0 in
$\PP^5$ is semistable with respect to the action of $SO(6)$. The
next proposition gives a criterion for an arbitrary quadric
$$
F = \sum_{i,j=1}^6 f_{ij}x_ix_j
$$
with $f_{ij}=f_{ji}$ for $i \neq j$ in $\PP^5$ to be semistable
with respect to the action of $SO(6)$.

\begin{prop} \label{prop3.2}
The quadric $F$ in $\PP^5$ is not semistable with respect to the
action of $SO(6)$ if and only it is equivalent under this action
to a quadric $Q = (q_{ij})$ with
$$
q_{ij} =0 \; for \; all \; 1 \leq i,j \leq 3 \quad and \quad
q_{14} = q_{15} = q_{16} = q_{25} = q_{35} = q_{45} = 0.
$$
\end{prop}
In other words, a quadric $F$ in $\PP^5$ is semistable with
respect to the action of $SO(6)$ if and only if there is no $g \in
SO(6)$ such that
\begin{equation} \label{eq3.2}
g^tFg = \left(\begin{array}{cccccc}0&0&0&0&0&0\\
                                          0&0&0&*&0&0\\
                                            0&0&0&*&*&0\\
                                             0&*&*&*&*&*\\
                                             0&0&*&*&*&*\\
                                             0&0&0&*&*&*
\end{array}\right).
\end{equation}
\begin{proof}
Choose the coordinates of $\PP^5$ in such a way that
\begin{equation} \label{eq3.3}
 G = \left(\begin{array}{cc}0&1_3\\
                            1_3&0
\end{array}\right).
\end{equation}
In classical terminology these coordinates are called {\it
Pl\"ucker coordinates}. Then the matrices
$diag(x_1,x_2,x_3,x_1^{-1},x_2^{-1},x_3^{-1})$ with $x_i \in \CC^*$
form a maximal torus of $SO(6)$. This fact and the form of the
Weyl group of $SO(6)$ imply that every 1-parameter subgroup of
$SO(6)$ is conjugate to one of the form
$$
\lambda: \CC^* \ra SO(6), \qquad t \mapsto
diag(t^{r_1},t^{r_2},t^{r_3},t^{r_4},t^{r_5},t^{r_6})
$$
with integers $r_1 \geq r_2 \geq r_3 \geq 0$ and $r_4 = -r_1,
r_5=-r_2, r_6=-r_3$ acting on the space of quadrics in the usual
way. In particular it acts on a monomial $x_ix_j$ of degree 2 by
$$
\lambda(t)(x_ix_j) = t^{-r_i-r_j}x_ix_j  \qquad \mbox{for} \qquad
1 \leq i \leq j \leq 6.
$$
Defining
$$
\mu(F,\lambda) = \max\{r_i + r_j \;| \; f_{ij} \ne 0\},
$$
the Hilbert-Mumford criterion implies that it suffices to show
that for a given quadric $Q = (q_{ij})$ there exists a $\lambda$
as above with $\mu(Q,\lambda) < 0$ if and only if the coefficients in the statement
of the proposition vanish.

It is easy to see that if there exists a 1-parameter group
$\lambda$ as above with $\mu(Q,\lambda) < 0$, then the
coefficients vanish. For example, if $q_{15} \neq 0$, then
$\mu(Q,\lambda) \geq r_1 + r_5 \geq r_2 - r_2 = 0$, a
contradiction.

Conversely, suppose all these coefficients vanish. So
$$
Q = 2q_{24}x_2x_4 + 2q_{34}x_3x_4 + 2q_{35}x_3x_5 + \sum_{i,j=4}^6
q_{ij}x_ix_j.
$$
Taking $r_1=3, r_2=2$ and $r_3 = 1$ we get
$$
\mu(Q,\lambda) =
max(2r_4,2r_5,2r_6,r_4+r_5,r_4+r_6,r_5+r_6,r_2+r_4,r_3+r_4,r_3+r_4)=
-1 < 0.
$$
This completes the proof of the proposition.
\end{proof}

\begin{cor} \label{cor3.3}
A quadratic complex $X$ is semistable with respect to the action
of $SO(G)$ if the discriminant $\Delta(X)$ has at least two
different roots, i.e. if its Segre symbol consists of at least 2
brackets.
\end{cor}
\begin{proof}
According to Proposition \ref{prop3.1} a quadratic complex in
Klein coordinates is semistable if and only if the quadric
$\Phi(X)$ of trace 0 is semistable. Changing to Pl\"ucker
coordinates the statement remains true, since according to
(\ref{eq2.1}) the matrix $A$ of the coordinate change satisfies
$A^t=A^{-1}$. Then $G$ is given by (\ref{eq3.3}). Let $F_0$ denote
the matrix $\Phi(X)$ transformed into Pl\"ucker coordinates, i.e.
$F_0 = A^t\Phi(X)A$.

According to Propositions \ref{prop3.1} and \ref{prop3.2}, $X$ is
not semistable if and only if $F_0$ is equivalent under the action
of $SO(6)$ to a matrix of the form of the right hand side of
(\ref{eq3.2}). Since the multiplicities of the roots of
$\Delta(X)$ stay the same under a change of coordinates, we may
even assume that $F_0$ is of this form. But then clearly
$\Delta(X)(\lambda,\mu) = \lambda^6$. In particular $\Delta$ has
only one root.
\end{proof}

\begin{rem}
{\em One could even work out which irreducible
quadratic complexes (see Lemma \ref{lem4.1} below) are not semistable. They are exactly those with Segre symbols 
$[6], [(51)], [(42)], [(33)], [(411)],$ $[(321)]$ and  $[(222)]$.
}
\end{rem}

\section{Moduli spaces of quadratic complexes}

In this section we construct the moduli spaces of quadratic
complexes with a fixed Segre symbol. First we need some
preliminaries. Recall that a quadratic complex is called {\it
irreducible} (respectively {\it non-reduced}), if it is
irreducible (respectively non-reduced) as a variety in $\PP^5$.

\begin{lem} \label{lem4.1}
(1) A quadratic complex is non-reduced if and only if its
Segre symbol contains a bracket of length 5;\\
(2) A quadratic complex is reducible if and only if its
Segre symbol contains a bracket of length $4$.\\
\end{lem}

\begin{proof}
(1): Suppose the Segre symbol of a quadratic complex $X$ contains
a bracket of length 5. It corresponds to a 4-cone, which is a
double plane $F$ in $\PP^5$. Hence $X = F \cap G$ is non-reduced.
Conversely, If $X$ is non-reduced, it follows from the Jacobi
criterion that the associated pencil contains a double plane. The
bracket corresponding to it is of length 5.

(2): Suppose the Segre symbol of $X$ contains a bracket of length
4. It corresponds to a 3-cone $F$. The directrix of $F$ is a
nonsingular quadric in $\PP^1$, i.e. consists of 2 different
points. Hence $F$ and thus $X$ is reducible. 

Conversely, suppose  $X$
reducible. No component of $X$ can be of degree 1, since otherwise the smooth quadric $G$ 
would contain a projective space of dimension 3. Hence $X = X_1 \cup X_2$ with smooth 3-dimensional quadrics $X_1$ and 
$X_2$. Let $H_i \;(\simeq \PP^4)$ denote the linear span of $X_i$ for $i=1$ and 2. Then $F_0 := H_1 \cup H_2$ is a quadric in $\PP^5$ 
such that
$$
X = F_0 \cap G.
$$
Two different $\PP^4$'s in $\PP^5$ intersect in a $\PP^3$. This means that $F_0$ is a 3-cone. The bracket 
corresponding to it is of length 4.  
\end{proof}

\begin{lem} \label{lem4.2}
Let $X$ and $X'$ be quadratic complexes with associated pencils
$\{\lambda F + \mu G \;|\; (\lambda: \mu) \in \PP^1\}$ and
$\{\lambda' F' + \mu' G \;|\; (\lambda':\mu') \in \PP^1\}$ and
roots $(\lambda_i:\mu_i)$ and $(\lambda'_j:\mu'_j)$ of the
corresponding discriminants.

Then the quadratic complexes $X$ and $X'$ are isomorphic if and
only if they have the same Segre symbol and there is an
automorphism of $\PP^1$ fixing $(0:1)$ and carrying
$(\lambda_i:\mu_i)$ to $(\lambda'_i:\mu'_i)$ for all $i$, where 
$(\lambda_i:\mu_i)$ and $(\lambda'_i:\mu'_i)$ correspond to brackets 
of the same type.
\end{lem}

\begin{proof}
For the proof we choose Segre coordinates. As quoted already in
Section 2, the corresponding pencils are isomorphic if and only if
they have the same Segre symbol and there is an automorphism
$(x_0:x_1) \mapsto (ax_0 + bx_1:cx_0+dx_1)$ of $\PP^1$ carrying
$(\lambda_i:\mu_i)$ into $(\lambda'_i:\mu'_i)$ for all $i$, where the brackets of 
$(\lambda_i:\mu_i)$ and $(\lambda'_i:\mu'_i)$ are of the same type. Now an
isomorphism of quadratic complexes maps $G$ onto $G$. This means
just that the automorphism of the associated pencils fixes the
point $(0:1)$ of $\PP^1$.
\end{proof}

\begin{cor} \label{cor4.3}
If $\sigma$ is a Segre symbol with at most 2 brackets, then all
quadratic complexes with Segre symbol $\sigma$ are isomorphic.
\end{cor}

\begin{proof}
Let $X_1$ and $X_2$ be quadratic complexes with Segre symbol
$\sigma$. Since the discriminants $\Delta(X_1)$ and $\Delta(X_2)$
admit $r \leq 2$ roots, there is an automorphism of $\PP^1$ fixing
$(0:1)$ and carrying the roots of $\Delta(X_1)$ onto the roots of
$\Delta(X_2)$. So the assertion follows from Lemma \ref{lem4.2}.
\end{proof}

Hence the moduli space of quadratic line complexes with a fixed
Segre symbol consists of a point only whenever the corresponding
discriminant admits at most 2 different roots. In particular the
varieties of cosingular complexes are not interesting for these Segre symbols. We
assume in the sequel that $\sigma$ is a Segre symbol with the
following 2 properties:
\begin{equation} \label{eq4}
\sigma \; \mbox{does not contain any brackets of length} \geq 4;
\end{equation}
\begin{equation}\label{eq5}
\sigma \; \mbox{consists of at least 3 brackets}.
\end{equation}
According to Lemma \ref{lem4.1}, (\ref{eq4}) implies that every
quadratic complex with Segre symbol $\sigma$ is irreducible and
reduced and (\ref{eq5}) means, as we shall see, that the
corresponding moduli space is positive dimensional. There are
exactly 23 Segre symbols with the properties (\ref{eq4}) and
(\ref{eq5}), see table 7.3 below. Let $\sigma$
be one of them. We want to construct the moduli space ${\cal
M}(\sigma)$ of quadratic complexes with Segre symbol $\sigma$.

\begin{lem} \label{lem4.4}
The quadratic line complexes with Segre symbol $\sigma$ are
parametrized by a quasiprojective subvariety $R(\sigma)$ of the
variety $LC \simeq \PP^{19}$ of all quadratic complexes.
\end{lem}

\begin{proof}
Clearly the quadratic complexes of Segre normal form with Segre
symbol $\sigma$ are parametrized by a quasiprojective variety
$\tilde{R}(\sigma)$. In fact, $\tilde{R}(\sigma) \simeq (\PP^1)^r
\setminus \{diagonals\}$, where $r$ denotes the number of brackets in
$\sigma$. Since every quadratic complex is isomorphic to one in
Segre normal form, $R(\sigma)$ is the image of the map
$$
SO(G) \times {\tilde R}(\sigma) \ra LC, \qquad (g,X) \mapsto g^tXg
$$
and as such a quasiprojective subvariety of $LC$.
\end{proof}
The group $SO(G)$ acts on the variety $R(\sigma)$ in an obvious
way. We have to determine the stabilizer of any $X \in R(\sigma)$.
\begin{lem} \label{lem4.5}
Let the Segre symbol $\sigma$ satisfy (\ref{eq4}) and (\ref{eq5})
and suppose that $\sigma$ consists of $r_i$ brackets of length $i$
for $i = 1,2$ and 3. Then the stabilizer of any quadratic complex
$X \in R(\sigma)$ in $SO(G)$ is of dimension
$$
\dim Stab(X) = r_2 + 3r_3,
$$
except in the case $\sigma = [(2,2),1,1]$, where it has dimension
2. In particular the dimension of the stabilizer depends only on
the Segre symbol $\sigma$ and not on the quadratic complex $X \in R(\sigma)$.
\end{lem}
\begin{proof}
A matrix $A \in SL(6)$ is in the stabilizer of $X$ if and only if
$A^tGA=G$ and $A^tFA=F.$ Since the pencil is in Segre's normal
form, $G=G^{-1}$, and these equations are equivalent to
\begin{equation} \label{eqn6}
 A^t=GA^{-1}G \quad \mbox{and} \quad GFA=AGF 
 \end{equation}
Suppose first that every cone is a zero cone. We have to show that the stabilizer is zero-dimensional.

It suffices to show that the stabilizer of every pair of blocks corresponding to a zero-cone of type $d$ is zero-dimensional.
The Segre normal form of these blocks $G_i$ of $G$ and $F_i$ of $F$ are given by the $d \times d$ matrices:
\[ G_i=\left(\begin{array}{ccccc}0&0&\dots&0&1\\
                                          0&0&\dots&1&0\\
                                            \dots&\dots&\dots&\dots&\dots\\
                                             0&1&\dots&0&0\\
                                             1&0&\dots&0&0
\end{array}\right)
\mbox{ and }\;
F_i=\left(\begin{array}{ccccc}0&0&\dots&1&\LL\\
                                        0&\dots&1& \LL&0\\
                                        \dots&\dots&\dots&\dots&\dots\\
                                        1&\LL&0&\dots&0\\
                                        \LL&0&0&\dots&0
\end{array}\right)\,
\]

Consider the $ d\times d$ matrix $A=\left(a_{ij}\right)$.
Computing $G_iF_iA$ and $AG_iF_i$, we see from the second equation of (\ref{eqn6}) that $A$ must be of the form:
\[ A=\left(\begin{array}{ccccc}a_{11}&0&\dots&0&0\\
                                          a_{21}&a_{11}&\dots&0&0\\
                                            \dots&\dots&\dots&\dots&\dots\\
                                            a_{(n-1),1}&a_{(n-2),1}&\dots&a_{11}&0\\
                                             a_{n1}&a_{(n-1),1}&\dots&a_{21}&a_{11}
\end{array}\right)\]
Using this, we deduce from $A^tG_iA=G_i$, that
\[G_i = \left(\begin{array}{cccccc}2a_{11}a_{n1}+2a_{21}a_{(n-1),1)}+\dots&\dots&\dots&\dots&2a_{11}a_{21}&a_{11}^2\\
                                          2a_{11}a_{(n-1),1}+2a_{21}a_{(n-2),1}&\dots&\dots&\dots&a_{11}^2&0\\
                                          \dots&\dots&\dots&\dots&\dots&\dots\\
                                           2a_{11}a_{31}+a_{21}^2&2a_{11}a_{21}&a_{11}^2&0&\dots&0\\
                                            2a_{11}a_{21}&a_{11}^2&0&0&\dots&0\\
                                             a_{11}^2&0&0&0&\dots&0
\end{array}\right)\]
This implies that $A = \pm id$ and thus the assertion in this case. Note that the result does not depend on $\lambda$, which
means that the dimension of the stabilizer does not depend on the chosen quadratic complex.  

Since the proof in the remaining cases is analogous, we omit the details. To be more precise, 
we distinguish the following cases and the proof always uses the equations (\ref{eqn6}).
 
First assume the pencil has $k_1$  $1$-cones none with a
bracket of type (2,2), and all other cones are $0$-cones. In this case
the stabilizer  is $k_1$-dimensional.

In case the pencil $\mathcal{P}$ has a $2$-cone, it can have
only this one $2$-cone according to our hypotheses. Segre's normal
form has therefore one block of either of these 2 forms:
$G_i= \mathds{1}$ and $F_i = \lambda \mathds{1}$ or
$G_i=\left(\begin{array}{ccc}0&1&0\\
                                1&0&0\\
                                 0&0&\mathds{1}
                                 \end{array}\right)\;
\mbox{ and }\;
 F_i=\left(\begin{array}{ccc}1&\LL&0\\
                            \LL&0&0\\
                                0&0&\LL\mathds{1}
                            \end{array}\right)$.
In the first case, the equation $A^t F_i A=F_i$ implies $A^t=A^{-1}$
and the stabilizer is $SO(3)$ which clearly 3-dimensional.
In the second  case the stabilizer also is seen to be 3-dimensional. It follows that if the pencil has $k_0$ $0$-cones and $k_1$
$1$-cones (none with a bracket  of type (2,2)) and $k_2 \; (= 1)$ 2-cones,
the stabilizer has dimension $k_1+ 3k_2$.

Finally consider the exceptional case $\sigma = [(2,2),1,1]$. The reason for the increase of the dimension of the stabilizer in this case
is that there are two $2 \time 2$ blocks corresponding to 
the same root of the discriminant and both this blocks are not diagonal. 
\end{proof}

As a consequence of Lemmas \ref{lem4.2}, \ref{lem4.4} and
\ref{lem4.5} we obtain
\begin{cor} \label{cor4.6}
With the assumptions of above, let $r = r_1 + r_2 + r_3$ be the
number of brackets in $\sigma$. Then we have
$$
\dim R(\sigma) = r + 13 - \dim Stab(X)
$$
where $X$ is any quadratic complex with Segre symbol $\sigma$.
\end{cor}
The main result of this section is the following theorem
\begin{teo} \label{teo4.7}
Let $\sigma$ be a Segre symbol satisfying (\ref{eq4}) and
(\ref{eq5}) consisting of $r$ brackets. Then the moduli space
${\cal M}(\sigma)$ of quadratic complexes exists and is a
quasiprojective variety of dimension $r-2$.
\end{teo}

\begin{proof}
According to Lemma \ref{lem4.4} the variety $R(\sigma)$
parametrizes all quadratic complexes with Segre symbol $\sigma$.
The group $SO(G)$ acts on $R(\sigma)$ and two quadratic complexes
are isomorphic if and only if they differ by this action.
According to Corollary \ref{cor3.3} every element of $R(\sigma)$
is semistable with respect to the action of $SO(G)$. Hence
according to \cite[Theorem 3.14]{Ne} a good
quotient ${\cal M}(\sigma)$ of $R(\sigma)$ modulo the action of
$SO(G)$ exists, is a quasiprojective variety and parametrizes the
closed orbits. Since by Lemma \ref{lem4.5} all orbits are of the
same dimension, all orbits are closed in $R(\sigma)$. Hence ${\cal
M}(\sigma)$ parametrizes the isomorphism classes of quadratic
complexes with Segre symbol $\sigma$. For the dimension we have
according to Corollary \ref{cor4.6}
$$
\dim {\cal M}(\sigma) = \dim R(\sigma) - \dim SO(G) + \dim Stab(X)
= r-2.
$$
\end{proof}

\section{The singular surface of a quadratic complex}

Recall that a point $x \in G$ represents a line in $\PP^3$ which
we denote by $l_x$. The points of $G$ which correspond to the
lines in $\PP^3$ passing through a particular point $p \in \PP^3$
form a plane $\alpha(p)$ contained in $G$. Similarly, the points
of $G$ corresponding to lines in $\PP^3$ lying in a plane $h
\subset \PP^3$ form a plane $\beta(h)$ contained in $G$. Therefore
we have two systems of planes on $G$ and in fact these are the
only planes in $G$. Two distinct planes of the same system meet
exactly in one point, while two planes of different systems are
either disjoint or meet exactly in a line. Conversely, every line
on $G$ is contained in exactly one plane of each of the two
systems. Following Newstead \cite{N}, we call the planes
$\alpha(p)$ {\it $\alpha$-planes} and the planes $\beta(h)$ {\it
$\beta$-planes}.

Now let $\sigma$ denote a Segre symbol satisfying properties
(\ref{eq4}) and (\ref{eq5}) and consider a quadratic line complex
$$
X = F \cap G,
$$
with Segre symbol $\sigma$. For a general point $p \in \PP^3$ the
intersection $\alpha(p) \cap F$ is a smooth conic in $\alpha(p)$.
The set
$$
S = \{p \in \PP^3 : {\mbox rk} (\alpha(p) \cap F) \leq 2\}
$$
is a surface in $\PP^3$, not necessarily irreducible. It is called
the {\it singular surface in $\PP^3$} associated to the complex
$X$. Clearly it does not depend on the quadric $F$ defining the quadratic complex, but only on the complex $X$ itself. The set
$$
R = \{p \in \PP^3 : {\mbox rk} (\alpha(p) \cap F) \leq 1\}
$$
is an algebraic subset of dimension $\leq 1$ of $S$. For $p \in S
\setminus R$, the intersection $\alpha(p) \cap F$ parametrizes two
different pencils of lines in $\PP^3$ intersecting in the common
point $p$, denoted as the {\it focus} of the two {\it cofocal
pencils}. For $p \in R$, either the intersection $\alpha(p) \cap
F$ parametrizes one pencil of lines in $\PP^3$ counted twice, or
the plane $\alpha(p)$ is contained in $X$.
In particular $S$ is the set of foci of pencils of lines in the complex $X$.\\

We now define a surface $\Sigma \subset X \subset \PP^5$ closely
related to $S$. For any $x \in X$ the line $l_x \in \PP^3$ is
called a {\it singular line} of the complex $X$ {\it at a point}
$p \in l_x$, if the plane $\alpha(p)$ is contained in the tangent
space $T_xF$. This means that the line $l_x$ belongs to more than
one pencil of the complex (or to one pencil counted twice).

If the line $l_x$ is singular at the point $p$, the point $p$ is
certainly contained in the surface $S$. Conversely, for $p \in S
\setminus R$ there is a unique singular line at $p$, namely the
line of intersection of the two corresponding cofocal pencils. If
$p \in R$, any line $l_x$ through $p$ is singular at $p$. We will
see that the set
$$
\Sigma := \{x \in X : l_x \; \mbox{is a singular line of the
complex}\;  X \}
$$
is a surface in $X$, not necessarily irreducible. We call it the
{\it singular surface in } $\PP^5$ associated to $X$. In order to
work out the relation between $\Sigma$ and $S$ we need the
following lemma.

\begin{lem} \label{lem5.1} Let $x \in \Sigma$.\\
{\em (a)} If $x$ is a smooth point of $X$, the line $l_x$ is
singular at exactly one point $p \in l_x$.\\
{\em (b)} If $x$ is a singular point of $X$, the line $l_x$ is
singular at any point $p \in l_x$.
\end{lem}

\begin{proof}
(a): Suppose the line $l_x$ is singular at the points $p \neq q$.
Since $\alpha(p) \cap \alpha(q) = x$, the linear span of
$\alpha(p)$ and $\alpha(q)$ in $\PP^5$ is the whole tangent space
$T_xG$. Hence $\alpha(p)$ and $\alpha(q)$ cannot be both contained
in $T_xF \neq T_xG$.\\
(b): If $x$ is a singular point of $X$, we have $T_xG \subset
T_xF$. So for any point $p \in l_x$, $\alpha(p) \subset T_xG
\subset T_xF$.
\end{proof}

\begin{rem} \label{rem5.2}
{\em If the $\alpha$-plane $\alpha(p)$ is contained in $X$, any
line passing through $p$ is singular exactly at $p$. If the
$\beta$-plane $\beta(h)$ is contained in $X$, any line in the plane
$h$ is singular at exactly one point $p$ and this gives a
bijection $\beta(h) \setminus Sing(X) \ra h \setminus Sing(S)$.}
\end{rem}

As a consequence of Lemma \ref{lem5.1} we can define a map
$$
\pi: \Sigma \setminus Sing(X) \ra S
$$
by associating to each $x \in \Sigma \setminus Sing(X)$ the unique
point $p \in l_x$ at which the line $l_x$ is singular.

\begin{lem} \label{lem5.3}
For any smooth point $x \in X$ the following conditions are
equivalent:\\
{\em (1)} $x \in \Sigma$;\\
{\em (2)} there is a $y \in G, \; y \neq x$ such that $T_xF =
T_{y}G$.
\end{lem}

\begin{proof}
(1) $\Rightarrow$ (2): Suppose $l_x$ is singular at $p \in \PP^3$,
i.e. $\alpha(p) \subset T_xF$. Since $x$ is a smooth point of $X$,
the tangent space $T_xF$ has projective dimension 4 and $G$
vanishes on $\alpha(p)$, which has projective dimension 2. This
implies that the restriction of $G$ to $T_xF$ is singular. So
$T_xF$ is tangent to $G$ at some point $y \in G$. Clearly $y \neq
x$, since otherwise $x$ would be a singular point of $X$.

(2) $\Rightarrow$ (1): Suppose $T_xF = T_yG$ for some $y \neq x$.
The line $\overline{xy} \subset \PP^5$ is contained in $G$ and
hence it is contained in exactly one plane of each system of
planes of $G$. Let $\alpha(p)$ the one corresponding to a point $p
\in \PP^3$. Then $\alpha(p) \subset T_xF$, i.e. $x \in \Sigma$.
\end{proof}

We may assume that $F$ is a smooth quadric, as is the case for the
Segre normal form. The following theorem shows that $\Sigma$ is a
complete intersection surface in $\PP^5$.

\begin{teo} \label{teo5.4}
Let $H$ denote the quadric defined by
the symmetric matrix $H = FG^{-1}F$. Then
\begin{equation}  \label{eq6}
\Sigma = F \cap G \cap H.
\end{equation}
\end{teo}
\begin{proof}
According to Lemma \ref{lem5.1} the line $l_x$ is singular for any
singular point $x$ of $X$, i.e. whenever $T_xF =T_xG$. Together
with Lemma \ref{lem5.3} we get that for any $x \in X$ we have: $x
\in \Sigma$ if and only if $T_xF = T_{x'}G$ for some point $x' \in
X$ (not necessarily different from $x$).

The dual coordinates of the tangent space $T_xF$, considered as a
point in ${\PP^5}^*$, are $x^* = Hx$, and this is tangent to $G$
if and only if $(x^*)^t G^{-1}x^* = 0$. This implies the
assertion.
\end{proof}

\begin{rem}
{\em According to the definition, $\Sigma$ depends only on $X$ and
not on the choice of $F$. This can be seen also from the
description in the theorem, since, if $F$ is replaced by $F +
\lambda G$, then $H$ is replaced by $H + 2\lambda F + \lambda^2
G$.}
\end{rem}

Theorem \ref{teo5.4} implies that the map $\pi : \Sigma \setminus
Sing(X) \ra S$ as defined above can be described in a way which
can be applied to compute the singular surface. For this we choose
the coordinates of $\PP^5$ in such a way that the matrix $G$
satisfies $G = G^{-1}$. This is the case for example for
Pl\"ucker-, Klein- and Segre coordinates. Then we have

\begin{prop} \label{prop5.5}
For any $x \in \Sigma \setminus Sing(X)$ the point $Fx$ is in $G$ and the map $\pi:
\Sigma \setminus Sing(X)  \ra S$ is given by
\begin{equation}  \label{eq7}
\pi(x) = l_x \cap l_{Fx}.
\end{equation}
\end{prop}

\begin{proof}
Suppose $x \in \Sigma$ is a smooth point of $X$.  Since $G =
G^{-1}$, we have according to (\ref{eq6}), $x \in \Sigma$ if and
only if
$$
x^tGx=0, \quad x^tFx = 0 \quad \mbox{and} \quad x^tFGFx = 0.
$$
The last equation can be read also as $(Fx)^tG(Fx) = 0$, which
means that $Fx \in G$.

Now with the coordinates $Y = (Y_1, \ldots, Y_6)$ of $\PP^5$,
$$
(Fx)^tY = 0
$$
is the equation of $T_xF$ as well as the equation of $T_{Fx}G$, so
that $T_xF = T_{Fx}G$. This implies that the line
$\overline{x,Fx}\subset \PP^5$ is contained in $G$. It follows
that the lines $l_x$ and $l_{Fx}$ intersect in a point $p \in
\PP^3$ and that the line $l_x$ is singular at $p$. This completes
the proof of the proposition.
\end{proof}

\section{Semistable quartics in $\PP^3$}

A quartic surface in $\PP^3$ is, up to a nonzero constant, defined
by a quartic form
\begin{equation} \label{eq6.1}
S = \sum_{i_0+i_1+i_2+i_3=4}
a_{i_0i_1i_2i_3}X_0^{i_0}X_1^{i_1}X_2^{i_2}X_3^{i_3}.
\end{equation}
As usual we denote the quartic surface and the quartic form by the
same letter. Two quartic surfaces $S$ and $S'$ are isomorphic if
there is a $\varphi \in SL(4)$ such that $\varphi(S) = S'$. This
defines an action of $SL(4)$ on the projective space $\PP^{34}$
parametrizing all quartic surfaces, which we want to analyze. The
line bundle $L = \cO_{\PP^3}(4)|S$ is $SL(4)$-linearizable, so
that we can speak about semistability of points in $\PP^{34}$.

\begin{lem} \label{lem6.1}
A quartic surface in $\PP^3$ is semistable with respect to the
action of $SL(4)$ if and only if it is not isomorphic to a
quartic (\ref{eq2.1}) with
$$
a_{4000} = a_{3100} = a_{3010}=a_{3001} = a_{2200} = a_{2110} =
a_{2101} = a_{2020} = a_{2011} = a_{2002} $$
$$
= a_{1300} = a_{1210} = a_{1201} = a_{1120} = a_{1111} = 0.
$$
\end{lem}

\begin{proof}
Consider the 1-parameter groups
$$
\lambda : \CC^* \ra SL(4), \qquad t \mapsto
diag(t^{r_0},t^{r_1},t^{r_2},t^{r_3})
$$
with integers $r_0 \geq r_1 \geq r_2 \geq r_3, \; \sum r_i =0$
acting on the space $\PP^{34}$ of quartics in $\PP^3$ in the usual
way. In particular it acts on a monomial
$x_0^{i_0}x_1^{i_1}x_2^{i_2}x_3^{i_2}$ of degree 4 by
$$
\lambda(t)(x_0^{i_0}x_1^{i_1}x_2^{i_2}x_3^{i_2}) =
t^{-(r_0i_0+r_1i_1+r_2i_2+r_3i_3)}x_0^{i_0}x_1^{i_1}x_2^{i_2}x_3^{i_2}
$$
Defining
$$ \mu(f,\lambda) = \max\{r_0i_0+r_1i_1+r_2i_2+r_3i_3 \;| \;
a_{i_0i_1 i_2 i_3} \ne 0\},
$$
the Hilbert-Mumford criterion implies that it suffices to show
that for a given quartic $f$ there exists a $\lambda$ as above
with $\mu(f,\lambda) < 0$ if and only if the coefficients
$a_{i_0i_1 i_2 i_3}$ of the proposition vanish.

It is easy to see that if there exists a one-parameter group
$\lambda$ as above with $\mu(f,\lambda) < 0$, then the
coefficients vanish. For example if $a_{3001} \neq 0$, then
$\mu(f,\lambda) \geq 3r_0+ r_3 = (r_0-r_1)+(r_0-r_2) \geq 0$.

Conversely, suppose all these coefficients vanish, and we take
$r_0=8,r_1=-1,r_2=-3,r_3=-4$, then $\mu(f,\lambda) = -1 < 0$. This
completes the proof of the lemma.
\end{proof}

\begin{prop}  \label{prop6.2}
A quartic surface $S \subset \PP^3$ is semistable with respect to
the action of $SL(4)$ if and only if it does not admit a triple
point whose tangent cone is a cone over a cuspidal plane cubic
(possibly degenerated).
\end{prop}

\begin{proof}
According to Lemma \ref{lem6.1}, $S$ is not semistable if and only
if it is isomorphic to a surface $S' \subset \PP^3$ with a triple
point at $e_0 = (1:0:0:0)$ whose tangent cone $TC(e_0)$ is given
in coordinates $y_i = \frac{x_i}{x_0}$ for $i=1,2$ and 3 by
\begin{equation} \label{eq6.1}
TC(e_0) = a_{1102}y_1y_3^2 + a_{1030}y_2^3 + a_{1021}y_2^2y_3 +
a_{1012}y_2y_3^2 + a_{1003}y_3^3.
\end{equation}
Considered as a plane projective curve, $TC(e_0)$ is a cubic with
a cusp at $(1:0:0)$ (or a degeneration of it). As an isomorphic
surface $S$ itself has a singularity of this type.

Conversely, since all cuspidal plane cubics are isomorphic, every
surface $S \subset \PP^3$ with a singularity of this type is
isomorphic to a surface $S \subset \PP^3$ with a triple point at
$(1:0:0:0)$ whose tangent cone is of the form (\ref{eq6.1}). By
Lemma \ref{lem6.1}, $S$ is not semistable.
\end{proof}

Using this we can construct various moduli spaces of quartic
surfaces in $\PP^3$. Recall that the projective space $\PP^{34}$
parametrizes all quartics in $\PP^3$. Equisingularity induces a
stratification of $\PP^{34}$ into locally closed algebraic
varieties (not necessarily irreducible). Certainly the group
action of $SL(4)$ restricts to an action on the strata. We may
assume that the dimension of the stabilizers is fixed on the
strata by refining the stratification if necessary. Since
semistability of a quartic depends only on the singularities of
the quartic, we may call a stratum {\it semistable} if one quartic
in it is.

Let $\sigma$ be a Segre symbol satisfying the assumptions
(\ref{eq4}) and (\ref{eq5}), i.e. $\sigma$ consists of at least 3
brackets and does not contain any bracket of length $\geq 4$. One deduces from the equations
of the normal forms given in \cite{J} and Section 7 that the quartics occurring as singular surfaces
of quadratic complexes with Segre symbol $\sigma$ all have the
same type of singularities. Let $Z_{\sigma} \subset \PP^{34}$
denote the corresponding stratum. We call the quartics of
$Z_{\sigma}$ {\it singular surfaces of type $\sigma$}.

\begin{lem} \label{lem6.3}
{\em (a)} The strata $Z_{\sigma}$ are locally closed subsets of the projective space $\PP^{34}$ 
parametrizing all quartics in $\PP^3$.

{\em (b)}
Any singular surface $S$ of a quadratic complex of type $\sigma$ 
is semistable with respect to the
action of $SL(4)$.
\end{lem}

\begin{proof}
(a): The singularities of the singular quartics of the quadratic complexes with Segre symbol $\sigma$ are (analytically) 
locally trivial in the sense of \cite{FK}. It is shown in \cite{FK} (see Corollary 0.2 and the proof of 0.3) that 
the locus of these quartics is analytically locally closed in the base space $\PP^{34}$. But $\PP^{34}$ being projective, this
implies the assertion.

(b): According to Proposition \ref{prop6.2}, $S$ is not semistable if
and only if it admits a triple point whose tangent cone is a cone
over a cuspidal cubic or a degeneration of it. Checking the equations of $S$ given in \cite{J} and Section 7 below, one sees
that this is not the case (it suffices to check this for the most
degenerate cases).
\end{proof}

\begin{teo} \label{teo6.4}
The moduli space of singular surfaces of type $\sigma$
$$
\cM_{ss}(\sigma) = Z_{\sigma}/SL(4)
$$
exists and is a quasiprojective variety.
\end{teo}

\begin{proof}
According to Lemma \ref{lem6.3} all elements of $Z_{\sigma}$ are
semistable. As for Lemma \ref{lem4.5} one checks that all stabilizers are of the same dimension. 
This implies that all orbits are of the same dimension. As in Theorem \ref{teo4.7} we conclude that
a geometric quotient $\cM_{ss}(\sigma) =
Z_{\sigma}/SL(4,\CC)$ exists and its points
parametrise the classes of isomorphic quartics in $Z_{\sigma}$.
\end{proof}

\begin{rem} \label{rem6.5}
{\em It is well known that the singular surface of a generic
quadratic complex is a {\it Kummer surface}, i.e. a quartic
surface in $\PP^3$, smooth apart from 16 ordinary double points.
Moreover every Kummer surface appears as the singular surface of a
generic quadratic complex. So in particular we constructed the
moduli space $\cM_{\kappa} := \cM_{ss}(\sigma)$ for $\sigma = [111111]$ 
of Kummer surfaces. Using the normal form for a Kummer
surface (see \cite{J} or \cite{Hu}) one checks that $\dim
Z_{\sigma} = 18$. On the other hand, it is easy to see that Kummer
surfaces have finite stabilizer in $SL(4)$ (see e.g.
\cite[Exercise V.5.1 (3)]{Be}). From this we conclude
$$
\dim \cM_{\kappa} = 3.
$$
Gonzalez-Dorrego \cite{GD} uses the normal form for Kummer
surfaces (see \cite[p.98]{J} or \cite[p.81]{Hu}) to construct the
moduli space $\cM_{\kappa}$ as follows: The normal forms
parametrize a 3-dimensional quasiprojective variety $\tilde{Z}$.
There is a finite group $N$ (an extension of the symmetric group
$S_6$ by ${\mathbb F}^4$), which is a subgroup of $SL(4)$ in a
natural way and thus acts on $\tilde{Z}$. The quotient
$\tilde{\cM}_{\kappa} = \tilde{Z}/N$ is the moduli space of Kummer
surfaces. The embedding $\tilde{Z} \ra Z_{\sigma}$ induces a
canonical isomorphism $\tilde{\cM}_{\kappa} \simeq \cM_{\kappa}$.}
\end{rem}

\section{The varieties of cosingular complexes}

Let $\sigma$ be a Segre symbol satisfying the assumptions
(\ref{eq4}) and (\ref{eq5}), i.e. $\sigma$ consists of at least 3
brackets and does not contain any bracket of length $\geq 4$. In
Theorem \ref{teo4.7} we constructed the moduli space $\cM_{qc}(\sigma)$
of quadratic complexes of type $\sigma$ and in Theorem
\ref{teo6.4} the moduli space $\cM_{ss}(\sigma)$ of quartic surfaces
of type $\sigma$.

As above let  $R(\sigma)$ and $Z_{\sigma}$ denote the spaces
parametrizing quadratic complexes and singular surfaces of type
$\sigma$ as in sections 4 and 6. In Section 5 we associated to
every quadratic complex in $R(\sigma)$ a singular surface in
$Z_{\sigma}$. This induces a map
$$ \pi: R(\sigma) \ra Z_{\sigma}. $$
which certainly is holomorphic. According to sections 4 and 6 the
groups $SO(6)$ and $SL(4)$ act on $R(\sigma)$ and $Z_{\sigma}$
in a natural way. Certainly these actions factorize via actions of
$PSO(6)$ and $PSL(4)$ respectively.

\begin{lem} \label{lem7.1}
There is an isomorphism $\iota: PSO(6) \ra PSL(4)$ such that
the map $\pi: R(\sigma) \ra Z_{\sigma}$ is equivariant with
respect to the actions of $PSO(6)$ and $PSL(4)$, i.e. 
$$
\pi(A
\cdot X) = \iota(A) \cdot \pi(X)
$$
for every $A \in PSO(6)$ and
$X \in R(\sigma)$.
\end{lem}

\begin{proof}
This is a well-known fact (see \cite[Section
19.1]{FH}).
In fact, the equivariance of the map $\pi$
can be used to the define the isomophism $\iota$: The points of
$\PP^3$ parametrize a family of planes, namely the
$\alpha$-planes, in the Pl\"ucker quadric $G$ and the action of
$PSO(6)$ on $\PP^5$ induces an action on this $\PP^3$. This
gives just the isomorphism $\iota$ of the lemma.
\end{proof}

\begin{rem}
{\em The equivariance of the map $\pi: R(\sigma) \ra Z_{\sigma}$
can also be expressed in terms of the actions of $SO(6)$ and
$SL(4)$: there is a surjective homomorphism $\kappa:
SL(4) \ra SO(6)$ with kernel of order 2 such that
$\pi(\kappa(\alpha) \cdot X) = \alpha \cdot \pi(X)$ for every
$\alpha \in SL(4)$ and $X \in R(\sigma)$.}
\end{rem}

Lemma \ref{lem7.1} implies that the map $\pi: R(\sigma) \ra
Z_{\sigma}$ induces a morphism of the corresponding moduli spaces,
which we denote by the same letter
$$ 
\pi: \cM_{qc}(\sigma) \ra \cM_{ss}(\sigma).
$$
Two quadratic complexes $X$ and $X'$ in $\cM_{qc}(\sigma)$ are called {\it cosingular} if their singular surfaces are isomorphic,
i.e. if $\pi(X) = \pi(X')$. The {\it variety $CS(X)$ of quadratic complexes cosingular to} $X$ is by definition the fibre of
the surface $\pi(X)$ under the map $\pi$:
$$
CS(X) := \pi^{-1}(\pi(X))
$$  
In the generic case $\sigma = [111111]$ the varieties $CS(X)$ have been investigated
by Klein in \cite{K}. It is the aim of this section to compute the
dimension of $CS(X)$ for a generic complex 
$X \in \cM_{ss}(\sigma)$ for every Segre symbol $\sigma$ satisfying equations (\ref{eq4}) and (\ref{eq5}).
The result is given in the table 7.3 below.\\

The order of the Segre symbols is chosen as in \cite[pp 230-232]{J}. We omit here, however, the quadratic
complexes with Segre symbol with less than 3 brackets.
We only give the equation of the singular surface, 
when the equation is either not given or incorrect (i.e. there is a typo) in \cite{J}. 
For the other cases we refer to the corresponding section of \cite{J}. 
\newpage
$$
{\bf Table \; 7.3}
$$
\begin{center}
\begin{tabular}{c|c|c|c|c|c}
& Segre symbol $\sigma$ & singular surface $S$& $\dim \cM_{qc}(\sigma)$  & $\dim \cM_{ss}(\sigma)$ & $\dim \pi^{-1}(S)$\\
\hline \hline 1 & $[111111]$ & see \cite[p. 98]{J} & 4 & 3 & 1 \\
\hline 2 & [21111]& see \cite[No 171]{J} & 3 & 2 & 1 \\
\hline 3 & [3111] & see \cite[No 180]{J} & 2& 1  & 1\\
\hline 4 & [411] &  see \cite[No 194]{J} & 1& 0  & 1\\
\hline 5 & [2211] & see \cite[No 186]{J} & 2& 1 & 1\\
\hline 6 & [321] & see \cite[No 198]{J} & 1 & 0 & 1\\
\hline 7 & [222] &  see \cite[No 203]{J}& 1 & 0 & 1\\
\hline 8 & [(11)1111] & see \cite[No 162]{J} &3 & 2 & 1\\
\hline 9 & [(11)211] & see \cite[No 173]{J} &2 & 1 & 1\\
\hline 10 & [(11)31] & see \cite[No 182]{J} & 1& 0 & 1\\
\hline 11 & [(11)22] & see Case 11 below& 1 &0 & 1\\
\hline 12 & [(21)111] & see \cite[No 172]{J} & 2& 1 & 1\\
\hline 13 & [(21)21] & see \cite[No 188]{J}& 1 & 0 & 1\\
\hline 14 & [(31)11] & see Case 14 below  &1 & 0 &1 \\
\hline 15 & [(22)11] & see \cite[No 187]{J} & 1& 0 & 1\\
\hline 16 & [(11)(11)11] & see \cite[No 167]{J} & 2& 1 & 1\\
\hline 17 & [(11)(11)2] & see \cite[No 177]{J} & 1& 0 & 1\\
\hline 18 & [(21)(11)1] & see Case 18 below & 1& 0 & 1\\
\hline 19 & [(11)(11)(11)] & coordinate tetrahedron &1 & 0 & 1\\
\hline 20 & [(111)111] & see Case 20 below & 2& 0 & 2\\
\hline 21 & [(111)(11)1] & same as in Case 19 &1 & 0 & 1 \\
\hline 22 & [(111)21] & see \cite[No 176]{J} &1 & 0 &1 \\
\hline 23 & [(211)11] & see Case 23 below &1 & 0 & 1\\

\end{tabular}
\end{center}

\noindent {\bf Case 11}: $\sigma = [(11)22]$\\
The equations for $G$ and $F$ are
$$
G = x_1^2 + x_2^2 + 2x_3x_4 + 2x_5x_6
$$
$$
F= \lambda_1(x_1^2 + x_2^2) + 2\lambda_2x_3x_4 + 2\lambda_3x_5x_6
+ x_3^2 + x_5^2
$$
and the equation of the singular surface is
$$
S = \lambda_1^2y_1^2y_3^2 + (\lambda_1-\lambda_2)^2y_1^2y_4^2
-4\lambda_1\lambda_2(\lambda_1-\lambda_2)y_1y_2y_3y_4.
$$

\noindent {\bf Case 14}: $\sigma = [(31)11]$\\
The equations of
$G$ and $F$ are
$$
G = x_1^2 + x_2^2 + x_3^2 + x_5^2 + 2x_4x_6
$$
$$
F = \lambda_1x_1^2 + \lambda_2x_2^2 + \lambda_3(x_3^2 + x_5^2 + 2x_4x_6) + 2x_4x_5.
$$
Then the equation of the singular surface $S$ is
$$
S= (\lambda_1 - \lambda_2)(y_1^4 + y_4^4) + 8(\lambda_1 - \lambda_3)(\lambda_2 - \lambda_3)y_1y_4(y_1y_3-y_2y_4) 
+ 2(\lambda_1 + \lambda_2 - 2 \lambda_3)y_1^2y_4^2.
$$


\noindent {\bf Case 18}: $\sigma = [(21)(11)1]$\\ The equations of
$G$ and $F$ are
\begin{equation} \label{eqn1}
G = x_1^2 + x_2^2 + x_3^2 + x_4^2 + 2x_5x_6
\end{equation}
$$
F = \lambda_1(x_1^2 + x_2^2) + \lambda_3x_3^2 + \lambda_4(x_4^2 +
2x_5x_6) +x_5^2.
$$
Then $S$ is given by
$$
S = (\lambda_1 - \lambda_4)(\lambda_3 -
\lambda_4)(y_1y_3+y_2y_4)^2 - 4(\lambda_3- \lambda_1)y_1^2 y_4^2.
$$

\noindent {\bf Case 20}: $\sigma = [(111)111]$\\ The quadric $G$
is given by $\sum_{i=1}^{6} x_i^2$ whereas
$$
F = \lambda_1(x_1^2 + x_2^2 + x_3^2)+ \lambda_4 x_4^2 + \lambda_5
x_5^2 + \lambda_6 x_6^2.
$$
Then the equation of $S$ is
$$
S  = (y_1y_3 - y_2y_4)^2.
$$

\noindent {\bf Case 23}: $\sigma = [(211)11]$\\
Let $G$ be as in (\ref{eqn1}) and $F$ be given by
$$
F = \lambda_1 x_1^2 + \lambda_2 x_2^2 + \lambda_3(x_3^2 + x_4^2 +
2x_5x_6) + x_5^2.
$$
Then
$$
S = y_1^2y_4^2.
$$

\begin{proof} For the proof of the second column, the equations of the singular surfaces, we applied the classical 
method as outlined in \cite{J}, but using Maple 9.5. The third column is a consequence of Theorem \ref{teo4.7}.
For the proof of the fourth column we showed in all the cases where $\dim \cM_{ss}(\sigma) =0$, again using Maple,  
that any 2 singular surfaces of that type are isomorphic. Since it is well known that the moduli space of Kummer 
surfaces is 3-dimensional, we can conclude the remaining dimensions by general arguments. To give an example, the variety
$\cM_{ss}([21111])$ is in the closure of the irreducible variety $\cM_{ss}([111111])$ and the variety $\cM_{ss}([(21)111])$
is in the closure of $\cM_{ss}([21111])$. So $\dim \cM_{ss}([(21)111]) \leq 1$. But $\cM_{ss}([(21)111])$ cannot be of dimension 0, 
since it is irreducible and the 0-dimensional variety $\cM_{ss}([(21)(11)1])$ is in its closure. Hence $\dim \cM_{ss}([(21)111]) = 1$
and we can conclude $\dim \cM_{ss}([21111]) = 2$. Finally the dimension $CS(X)$ for a general $X \in \cM_{qc}(\sigma)$ is given as the 
difference of the dimensions of $\cM_{qc}(\sigma)$ and $\cM_{ss}(\sigma)$.  
\end{proof}

\noindent
{\bf Remark 7.4.}
It is easy to work out the singularities of the singular surface in every case.
For example the singularities of the singular surface of type $[21111]$ consists of one line and 8 points. 
The classical authors in general did not mention the points (see \cite{J}, some singular points are, however, 
computed in \cite{W}).

\section{Cosingular complexes in the generic case}

In this section we investigate the variety of cosingular complexes of a
generic quadratic complex. Let $\sigma = [1 \ldots 1]$ which we assume in the whole section. 
In particular $\cM_{ss}(\sigma)$ is the moduli space of Kummer surfaces. Since
$\pi : \cM_{qc}(\sigma) \ra \cM_{ss}(\sigma)$ is surjective (see Remark
\ref{rem6.5}), $\dim \cM_{qc}(\sigma) = 4$ and $\dim \cM_{ss}(\sigma) = 3$,
a general fibre of $\pi$ is of dimension 1. We will see that all fibres 
are curves in this case. The main ideas of this section are due to Klein (see \cite{K}). 
We reformulate his arguments using our set-up.\\

Consider a fixed generic complex $X = F \cap G$ in Segre normal
form, (which in this case is the same as the Klein normal form) i.e.
\begin{equation}  \label{eq7.0}
G = \sum_{i=1}^6 x_i^2=0 \quad \mbox{and} \quad F = \sum_{i=1}^6
\lambda_i x_i^2 = 0
\end{equation}
with $\lambda_i \in \CC$ pairwise different. For any $\rho \in \CC,\; \rho
\neq \lambda_i$ for $i = 1, \ldots,6$ consider the quadric
$F_{\rho}$ with equation
\begin{equation}  \label{eq7.1}
F_{\rho} = \sum_{i=1}^6 \frac{x_i^2}{\lambda_i - \rho} = 0.
\end{equation}

\begin{lem} \label{lem7.3}
Let $\Sigma$ (respectively $\Sigma_{\rho}$) denote the singular
surface in $\PP^5$ of the complex $X = F \cap G$ (respectively
$X_{\rho} = F_{\rho} \cap G$). Then there is an automorphism
$\varphi$ of $\PP^5$ such that $\varphi(\Sigma) = \Sigma_{\rho}$.
\end{lem}

\begin{proof} According to Theorem \ref{teo5.4}, the surface $\Sigma$ is the complete intersection
in $\PP^5$ with equations
\begin{equation} \label{eq7.2}
\sum_{i=1}^6 x_i^2 = \sum_{i=1}^6 \lambda_ix_i^2 = \sum_{i=1}^6
\lambda_i^2 x_i^2 = 0.
\end{equation}
Similarly $\Sigma_{\rho}$ is given by the equations
\begin{equation}  \label{eq7.3}
\sum_{i=1}^6 x_i^2 = \sum_{i=1}^6 \frac{x_i^2}{\lambda_i - \rho} =
\sum_{i=1}^6 \frac{x_i^2}{(\lambda_i - \rho)^2} = 0.
\end{equation}
Let the automorphism $\varphi$ of $\PP^5$ with $\varphi(x) = y$ be
defined by
$$
x_i = \frac{y_i}{\lambda_i - \rho} \quad \mbox{for} \quad i = 1,
\ldots 6.
$$
If $x \in \Sigma$, then
$$
\begin{array}{c}
0 = \sum x_i^2 = \sum \frac{y_i^2}{(\lambda_i - \rho)^2},\\
0 = \sum \lambda_i x_i^2 = \sum \lambda_i \frac{y_i^2}{(\lambda_i
- \rho)^2} - \rho \sum \frac{y_i^2}{(\lambda_i - \rho)^2} =
\sum \frac{y_i^2}{\lambda_i - \rho},\\
0 = \sum \lambda_i^2 x_i^2 = \sum \lambda_i^2
\frac{y_i^2}{(\lambda_i - \rho)^2} 
- 2 \rho \sum \lambda_i \frac{y_i^2}{(\lambda_i - \rho)^2} 
+ \rho^2 \sum \frac{y_i^2}{(\lambda_i - \rho)^2} = \sum y_i^2.
\end{array}
$$
Hence $y = \varphi(x) \in \Sigma_{\rho}$. Similarly one checks
$\varphi^{-1}(\Sigma_{\rho}) \subset \Sigma$.
\end{proof}

Lemma \ref{lem7.3} implies that all quadratic complexes $X_{\rho}$
are contained in the fibre $\pi^{-1}S = \pi^{-1} \pi(X)$.

\begin{lem} \label{lem7.4}
The quadratic complex $X$ is the limit of the complexes $X_{\rho}$
as $\rho \ra \infty$.
\end{lem}

\begin{proof}
Fix $\rho_0 \in \CC \setminus \{\lambda_1, \ldots , \lambda_6\}$.
The complex $X_{\rho}$ can be described as $X_{\rho} = F'_{\rho}
\cap G$ with
$$
F'_{\rho} = \rho ( G + (\rho_0 + \rho) F_{\rho}) = \rho
\sum_{i=1}^6 \frac{\lambda_i + \rho_0}{\lambda_i - \rho} x_i^2 =
\sum_{i=1}^6 \frac{\lambda_i + \rho_0}{\frac{\lambda_i}{\rho} +
1}x_i^2.
$$
Hence $\lim_{\rho \ra \infty} F'_{\rho} = \sum_i \lambda_i x_i^2 + \rho_0
\sum_i x_i^2 = F + \rho_0 G$ which gives the assertion.
\end{proof}

For any $X \in \cM_{qc}(\sigma)$ consider the family of quadratic
complexes
$$
{\cal C}_X = \{ X, X_{\rho} \in \cM_{qc}(\sigma) \;|\;X_{\rho} = F_{\rho}
\cap G \; \mbox{with} \; F_{\rho} \; \mbox{as in} \; (\ref{eq7.1})
\}.
$$
Certainly the index $\rho$ in $X_{\rho}$ depends on the chosen
quadric $F_{\rho}$. However we have

\begin{prop} \label{prop7.5}
For a quadratic complex $X' = F' \cap G \in \cM_{qc}(\sigma)$ with $ F' = \sum_{i=1}^6
\lambda'_i x_i^2$ the following statements are equivalent\\
{\em (1)} $X' \in {\cal C}_X$;\\
{\em (2)} There is a matrix $A = (a_{ij}) \in SL(2)$ such that
for $i=1, \ldots,6$,
$$\lambda'_i = \frac{a_{11}\lambda_i + a_{12}}{a_{21}\lambda_i +
a_{22}}.$$
\end{prop}

\begin{proof}
(1) $\Rightarrow$ (2): If $X' = X$, choose $A = {\bf 1}_2$. So let
$X' = X_{\rho}$ with $\rho \neq \lambda_i$ for $i = 1, \ldots, 6$.
Then $A = \left( \begin{array}{cc} 0&-1\\
                               1&-\rho
                               \end{array} \right) \in SL(2)$
with $X' = F' \cap G$ where $F' = \sum_i \frac{-1}{\lambda_i -
\rho} x_i^2$, i.e. $\lambda'_i = \frac{-1}{\lambda_i - \rho}$.

(2) $\Rightarrow$ (1): Let $A \in SL(2)$ be as in (2). If
$a_{21} =0$, then $F' = \frac{a_{11}}{a_{22}}F +
\frac{a_{12}}{a_{22}}G$ and thus $X' = X$. If $a_{21} \neq 0$,
then
$$
\lambda'_i = \frac{\frac{a_{11}}{a_{21}}(a_{21} \lambda_i + a_{22})
+ a_{12} - \frac{a_{11}a_{22}}{a_{21}}}{a_{21}\lambda_i + a_{22}}
= \frac{a_{11}}{a_{21}} - \frac{1}{a_{21}^2(\lambda_i +
\frac{a_{22}}{a_{21}})}
$$
and thus $F' = -\frac{1}{a_{21}^2} F_{(-\frac{a_{22}}{a_{21}})} +
\frac{a_{11}}{a_{21}}G$, i.e. $X' \in {\cal C}_X$.
\end{proof}

Call two quadratic complexes $X'$ and $X$ {\it equivalent}, i.e.
$X' \sim X$, if and only if $X' \in {\cal C}_X$. Since $SL(2)$
is a group, we obtain as an immediate consequence of the
proposition

\begin{cor} \label{cor7.6}
$'' \sim ''$ is an equivalence relation on the set $\cM_{qc}(\sigma)$.
\end{cor}

Using this we can determine the fibres of the morphism $\pi:
\cM_{qc}(\sigma) \ra \cM_{ss}(\sigma)$.

\begin{teo} \label{teo7.7}
For any quadratic complex $X \in \cM_{qc}(\sigma)$,
$$
\pi^{-1}\pi(X) = {\cal C}_X.
$$
\end{teo}

\begin{proof}
As quotients of normal varieties by reductive groups, the moduli
spaces $\cM_{qc}(\sigma)$ and $\cM_{ss}(\sigma)$ are normal varieties. The generic
fibre of $\pi$ being irreducible, every fibre is connected by
Zariski's connectedness theorem. Now Corollary \ref{cor7.6}
implies the assertion.
\end{proof}

As a consequence of Proposition \ref{prop7.5} and Theorem
\ref{teo7.7} we obtain

\begin{cor}  \label{cor7.8}
Two quadratic complexes of $\cM_{qc}(\sigma)$ have the same singular
surface if and only if they are isomorphic as pencils of quadrics.
\end{cor}

\begin{prop}  \label{prop7.9}
Let $X  \in \cM_{qc}(\sigma)$ be a quadratic complex. For any $\rho_1,
\rho_2 \in \CC, \rho_1 \neq \rho_2$, the quadratic complexes
$X_{\rho_1}$ and $X_{\rho_2}$ are non-isomorphic and not
isomorphic to $X$.
\end{prop}

\begin{proof}
It suffices to show that for any $\rho \in \CC, \; \rho \neq \lambda_i$ for all $i$ the complex
$X_{\rho}$ is not isomorphic to $X$. Suppose $X = F \cap G$ with
$F$ and $G$ as in (\ref{eq7.0}). Then $X_{\rho} = F_{\rho} \cap G$
with $F_{\rho}$ as in (\ref{eq7.1}).

The roots of $\det(\lambda F + \mu G) = \prod_{i=1}^6 (\lambda
\lambda_i + \mu)$ are $(\lambda: \mu) = (1:-\lambda_i)$ and the
roots of $\det(\lambda F_{\rho} + \mu G) = \prod_{i=1}^6 (\lambda
\frac{1}{\lambda_i - \rho} + \mu)$ are $(\lambda: \mu) =
(1:-\frac{1}{\lambda_i - \rho})$. Hence according to Lemma
\ref{lem4.2} the complexes $X$ and $X_{\rho}$ are isomorphic if
and only if there is a permutation $\sigma$ of the indices $1,
\ldots, 6$ such that the system of linear equations in $c$ and $d$
\begin{equation}   \label{eq7.4}
c - d \lambda_i = - \frac{1}{\lambda_{\sigma(i)} - \rho}   \qquad
\mbox{for} \qquad i = 1, \ldots, 6
\end{equation}
admits a solution.\\
Since $\lambda_1 \neq \lambda_2$, there is a unique solution of
the first 2 equations. If this would be also a solution for $i =
3, \ldots, 6$, we would have
\begin{equation}  \label{eq7.5}
\det \left( \begin{array}{ccc} 1& -\lambda_1& -
\frac{1}{\lambda_{\sigma(1)} - \rho}\\
         1& -\lambda_2& - \frac{1}{\lambda_{\sigma(2)} - \rho}\\
         1& -\lambda_i& - \frac{1}{\lambda_{\sigma(i)} - \rho}\\
         \end{array} \right) = 0
         \end{equation}
for $i = 3, \ldots 6$. But for fixed $\lambda_1$ and $\lambda_2$
this has at most 2 solutions in $\lambda_i$ and
$\lambda_{\sigma(i)}$, whereas $\lambda_3, \ldots, \lambda_6$ are
4 different values. Hence for any permutation $\sigma$ the linear
system (\ref{eq7.4}) is unsolvable, which implies the assertion.
\end{proof}

\begin{center}

\vskip20pt

$\begin{array} {c}

\hbox{\footnotesize Departamento de Matem\'atica, UFMG}\\

\hbox{\footnotesize Belo Horizonte, MG 30161--970, Brasil.}

\\

\hbox{\footnotesize dan@mat.ufmg.br}

\end{array}$\quad\quad

$\begin{array} c

\hbox{\footnotesize Mathematisches Institut}\\

\hbox{\footnotesize Bismarckstr. $1\frac{1}{2},$ 91054, Erlangen}

\\

\hbox{\footnotesize lange@mi.uni-erlangen.de}

\end{array}$

\end{center}

\end{document}